\documentclass[twocolumn]{revtex4-2}
\usepackage{graphicx}
\usepackage{graphicx} 
\usepackage{amsmath}
\usepackage{amssymb}
\usepackage{color}
\usepackage{dcolumn}
\usepackage{bm}
\usepackage{mathptmx}
\usepackage{verbatim}
\usepackage{fdsymbol}
\usepackage{sidecap}
\usepackage{placeins}
\usepackage{xcolor}

\begin{document}

\title{On a stable torus in a 3D system with a saddle-focus.}
\author{Andrey L. Shilnikov}
\email{ashilnikov@gsu.edu}
\affiliation{Neuroscience Institute and Department of Mathematics \& Statistics, Georgia State University, \\
100 Piedmont Ave., Atlanta, GA 30303, USA.}
\author{Leonid P. Shilnikov}
~\\

\begin{abstract}
This paper proposes a conceptual model for the onset of a stable torus near a saddle-focus equilibrium. This bifurcation scenario is typical of slow-fast systems that generate elliptic bursting in a variety of neuronal models in mathematical neuroscience. Variants of the model also capture other dynamical regimes recurring in a neighborhood of the saddle-focus. We also discuss homoclinic bifurcations for which the model assumptions are feasible. 
\end{abstract}

\maketitle

\section{Preface}
\label{intro}

The text below is a draft of an unfinished paper on which my father, Leonid Pavlovich Shilnikov, worked in December 2010. Figure~\ref{LP} shows him at his desk while drafting this manuscript. He was never able to complete it; he passed away the following year, shortly after his seventy-seventh birthday. For many years after his death, I, Andrey, could not bring myself  to return to his handwritten notes and attempt to finish the work that he had begun. Now, nearly fourteen years later, with the much-needed help of Dima Turaev, I was finally able to resume and bring it to completion.
  
The idea of this paper was inspired by our conversations about research that my Ph.D. student, Jeremy Wojcik (now a well-established researcher), and I were conducting at that time. It concerns a three-dimensional slow-fast system in the context of so-called elliptic bursting, which frequently arises in phenomenological and biologically plausible models in mathematical neuroscience (see Refs.~\cite{WALS1,WALS2} and references therein). The bifurcation mechanism underlying an elliptic burster typically involves three key elements: a saddle-node bifurcation of two periodic orbits (one stable and one unstable), followed by a subcritical Andronov-Hopf bifurcation in the fast two-dimensional subsystem. In combination with a slow passage governed by a one-dimensional slow equation, which induces a delayed loss of stability, this configuration gives rise to endogenous bursting dynamics that may be regular (Lyapunov stable) or even chaotic.
Figure~\ref{fig2}A illustrates the geometry of the slow motion, or critical manifolds, in the three-dimensional phase space of an elliptic burster generated by a three-dimensional extension of the 2D canonical Van der Pol (VdP) relaxation oscillator. Note that the Van der Pol (VdP) relaxation oscillator -- given by the first two lines in Eqs.~(\ref{fhn}) -- is commonly referred to in the recent biosciences modeling literature as the FitzHugh-Nagumo model. We note that both independent authors, R.~FitzHugh and J.~Nagumo, explicitly acknowledged the original two-dimensional ODE system and referred to it as the VdP model~\cite{nagumo,fitzhugh}. We also note that the hysteresis underlying the well-known VdP relaxation oscillations is a pivotal feature of many biologically plausible Hodgkin-Huxley-type models~\cite{chaos, torus} in neuroscience. Needless to add, chaotically recurrent oscillations observed in neural systems are often linked to, and triggered by, the Shilnikov saddle-focus, which is also discussed below.

The simplest toy model describing elliptic bursting (its ``voltage'' trace is shown in Fig.~\ref{fig2}B) is given by the following ODE system with a single nonlinear term:
 \begin{figure}[t!]
\begin{center} 
\includegraphics[width=0.48\textwidth]{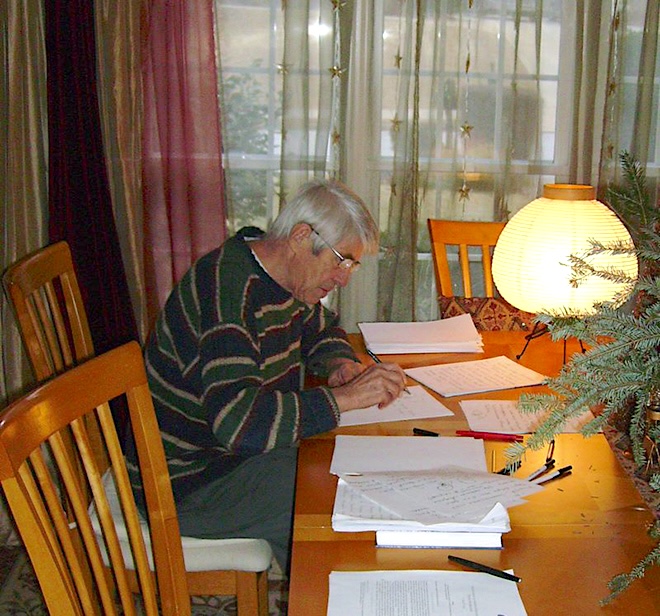}
\end{center}
\caption{Leonid Shilnikov drafting this paper at Andrey's home shortly before Christmas in 2010.}  
\label{LP}
\end{figure}
 \begin{equation} \label{fhn}
 \begin{array}{rclcl}
 v^\prime &= & v-v^3/3-w+y+I,\\
 w^\prime &= & \delta(0.7+v-0.8w),\\
 y^\prime &= & \mu(c-y-v),
 \end{array}
\end{equation}
and with several positive constants, where a small one, $\mu$, determines the pace of the slow $y$-variable. The slow variable $y$ becomes frozen in the singular limit $\mu=0$. Here we employ $c$ as the primary bifurcation parameter of the model, variations of which elevate or lower the slow nullcline $y^\prime=0$, which is the (grey) plane $v=y-c$ in the 3D phase space in Fig.~\ref{fig2}A. We have $y'>0$ below this plane and $y'<0$ above. The two fast equations in Eqs.~(\ref{fhn}) describe a relaxation oscillator in a plane, provided $\delta$ is small. The periodic oscillations in the fast subsystem are due to hysteresis induced by the cubic nonlinearity in the first ``voltage'' equation of the model. It is well known that oscillations cease with the emergence of a stable equilibrium state on either the upper or lower (but not the middle) branch of the cubic nullcline. We note that the factor 1/3 is crucial for the fast subsystem to demonstrate a sub-critical AH bifurcation; without it, the bifurcation becomes super-critical. This indicates that the first Lyapunov coefficient, whose sign determines the type of criticality, is close to zero and may change with parameter variations of the full, coupled system. In the sub-critical case, a stable equilibrium in the fast subsystem loses its stability after a repelling limit cycle (LC) collapses into it. Both,  repelling and stable, LCs emerge through a saddle-node bifurcation. This bifurcation occurs in the fold of the slow-motion manifold M$_{\rm LC}$ in the phase space of the full system. The inner cone is foliated by repelling LCs in the fast subsystem, which should be treated as saddle periodic orbits (PO) in the 3D phase space of the whole model. The outer surface of M$_{\rm LC}$ is comprised of stable POs. 

 \begin{figure}[t!]
\begin{center} 
\includegraphics[width=0.35\textwidth]{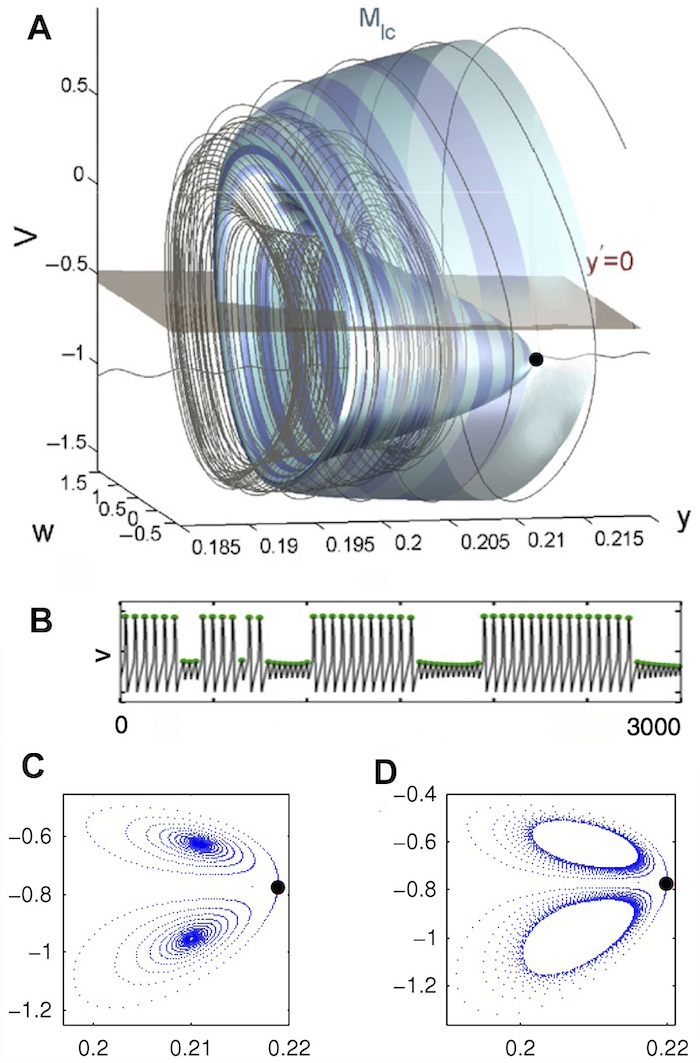}
\end{center}
\caption{(A) Geometry of an elliptic burster in phase space. The surface M$_{\rm LC}$ is the slow-motion manifold composed of stable (outer part) and unstable (inner cone) periodic orbits of the fast subsystem of the model~(\ref{fhn}). Its fold corresponds to a saddle-node bifurcation of periodic orbits. The black dot marks a subcritical AH bifurcation, to the left and right of which the quiescent equilibrium is a stable and unstable focus, respectively. Also shown is a trajectory converging to a stable two-dimensional torus in the three-dimensional phase space.
(B) Voltage waveform of the elliptic burster. (C,D) Trajectories spiraling toward a stable equilibrium and a smooth invariant circle on a cross-section transverse to M$_{\rm LC}$, corresponding to a stable periodic orbit and a two-dimensional torus, respectively.}  
\label{fig2}
\end{figure}
The full system can exhibit bursting, tonic spiking, or quiescent behavior for different values of the parameter $c$. Vertical shifts of the slow nullcline control transitions between these activity types by changing how and where it intersects the slow-motion manifold in phase space.
Loosely speaking, when the grey slow nullcline intersects the manifold M$_{\rm LC}$ through its upper branch, the model produces tonic spiking. When it intersects the unstable inner cone, the dynamics eventually becomes that of an elliptic burster. As the nullcline moves downward through the fold of M$_{\rm LC}$, the model may display either a period-doubling cascade initiated by the loss of stability of a periodic orbit, or beating and quasi-periodic oscillations associated with the emergence of a stable two-dimensional torus in phase space; see Figs.~\ref{fig2}C and D. These panels show a vertical cross-section transverse to the slow-motion manifold M$_{\rm LC}$ that is traversed by trajectories of the system.
In Fig.~\ref{fig2}C, trajectories initiated near the tip of the cone spiral toward two fixed points on the section, which represent the images of a single stable periodic orbit located inside the cone. In Fig.~\ref{fig2}D, this periodic orbit has become unstable, and trajectories instead converge alternately to two smooth, closed invariant circles, clearly indicating the presence of a stable two-dimensional torus in the vicinity of the cone in the phase space of the model.
Note also that when the slow nullcline (shown as the grey plane in the figure) passes exactly through the tip of the inner cone, the unique equilibrium of the full system has characteristic exponents $(\pm i\omega,0)$ in the limit $\mu=0$. Moreover, when the nullcline intersects the equilibrium manifold of the fast subsystem to the left or right of this tipping point, the equilibrium of the full system is, respectively, a stable focus or a saddle-focus. Further details and additional references can be found in Ref.~\cite{torus}, which describes how torus bifurcations enrich both simple and complex dynamics in a variety of neuronal models.

\begin{figure}[t!]
\begin{center} 
\includegraphics[width=0.48\textwidth]{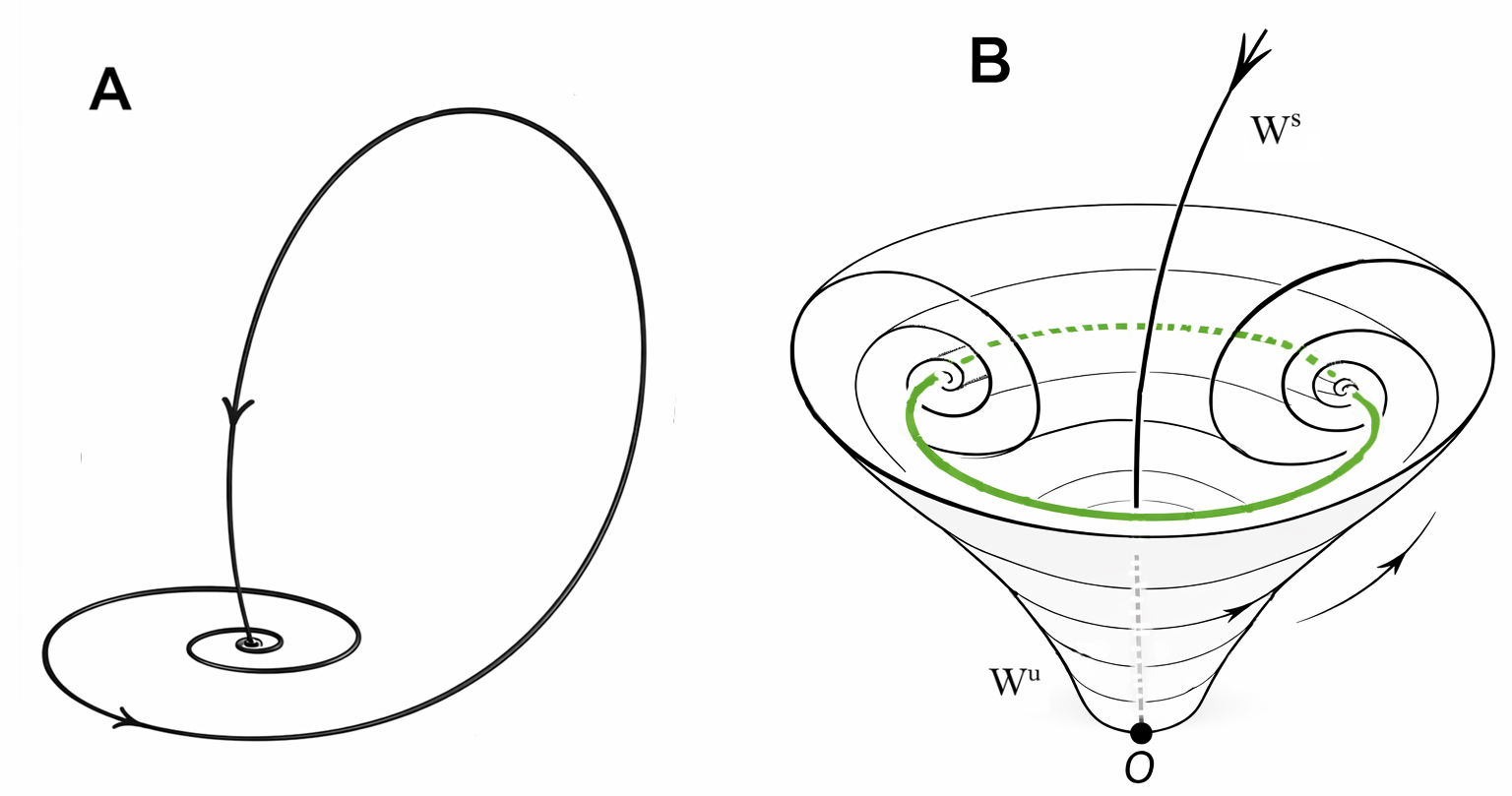}
\end{center}
\caption{(A) A homoclinic loop to a saddle-focus with a two-dimensional unstable manifold $W^{u}(O)$ and a one-dimensional stable manifold $W^{s}(O)$. in three-dimensional phase space. (B)  A sketch of a Shilnikov funnel generated by $W^{u}(O)$  of a saddle–focus: it widens by spiraling outward, forming a funnel-shaped surface, wrapping around onto an attractor: simple with regular or non-trivial with chaotic dynamics; here it is a stable periodic orbit (a green circle) -- just like as a period-doubling cascade had started.}  
\label{fig3}
\end{figure}

Here we arrive at the central motivation of the paper: how can a stable torus emerge in phase space near a homoclinic orbit to a saddle-focus (see Fig.~\ref{fig3}A)? Leonid Shilnikov’s long-term goal was to place this question within a broader conceptual paradigm aimed at describing canonical dynamical models whose behavior is confined to a saddle-focus funnel (or whirlpool) — an attracting, forward-invariant region of phase space organized and shaped by the two-dimensional unstable manifold of a saddle-focus, as sketched in Fig.~\ref{fig3}B. He outlined this program in the draft reproduced below and restricted attention to the three-dimensional case for the sake of simplicity. The main results were intended to be formulated as three propositions giving conditions for the existence of an invariant torus (Proposition I below) or, alternatively, nontrivial hyperbolic sets (Proposition II) or a stable periodic orbit (Proposition III) within the funnel. Proofs of Propositions I and II can be obtained by applying the Afraimovich-Shilnikov theory of torus breakdown~\cite{ALP3,ALP1,ALP2,torusbd,AS:1983,ALP4}, which is one of Leonid Shilnikov’s ground-breaking contributions, alongside the theory of saddle-focus homoclinic loops~\cite{Sciheritage,LPmathDC,LPbook17,book}. Proposition III is similar in spirit to the blue-sky catastrophe~\cite{blue,blue1,Mmo2005,blue2,blue3}.

The funnel construction was proposed by Leonid Shilnikov in~\cite{LP1986} as a geometric model for so-called spiral chaotic attractors \cite{AS:1983,LP1997,sfbif,Sciheritage}. The funnel can be viewed as a counterpart to the classical Afraimovich-Bykov-Shilnikov geometric model of the Lorenz attractor~\cite{ABS:1977,ABS:1983}. Leonid’s vision was that transitions from stationary to turbulent regimes in hydrodynamics, as well as in other classes of systems, are naturally associated with a saddle-focus that emerges after an Andronov-Hopf bifurcation and subsequently becomes part of the attractor. The saddle-focus' incorporation leads to characteristic spiraling episodes in dynamics. From this perspective, the basic unifying features of the transitions to turbulence can be captured by the Shilnikov's funnel framework.

As Leonid proposed in~\cite{LP1986}, the dynamics inside the funnel can be described by certain classes of annulus maps. The text below introduces the main types of such maps and their properties in a particular singular limit, and outlines several case studies in which this theory may be applicable. Since the original draft is sketchy in places, we have added brief comments, indicated in square brackets, where clarification was needed.
 
\section{Original draft}\label{section_model} 

We will consider a three-dimensional system
\begin{equation}\label{eq1}
\begin{aligned}
\dot x & = \rho x - \omega y+ \cdots,\\
\dot y & = \omega x + \rho y +  \cdots,\\
\dot z &= -\lambda z + \cdots, 
\end{aligned}
\end{equation}
where $\rho>0$, $\lambda >0$, $\omega \ge 0$, [the dots stand for the nonlinear terms].  Under broad conditions, the system near the origin $O$ can be reduced to a linear form
\begin{equation}\label{eq2}
\begin{aligned}
\dot x & = \rho x - \omega y,\\
\dot y & = \omega x + \rho y,\\
\dot z &= -\lambda z + \cdots, 
\end{aligned}
\end{equation}
or, in polar coordinates, 
\begin{equation}
\dot r  = \rho r , \qquad \dot \varphi = \omega, \qquad \dot z= -\lambda z. 
\label{eq3}
\end{equation}
[The coordinate origin in this linear system is a saddle-focus equilibrium state: it has a 1D stable manifold $W^s_O$ -- the $z$-axis ($x=y=0$), and a 2D unstable manifold $W^u_O$, which is the plane  $z=0$.] 

\begin{figure}[t!]
\begin{center} 
\includegraphics[width=0.499\textwidth]{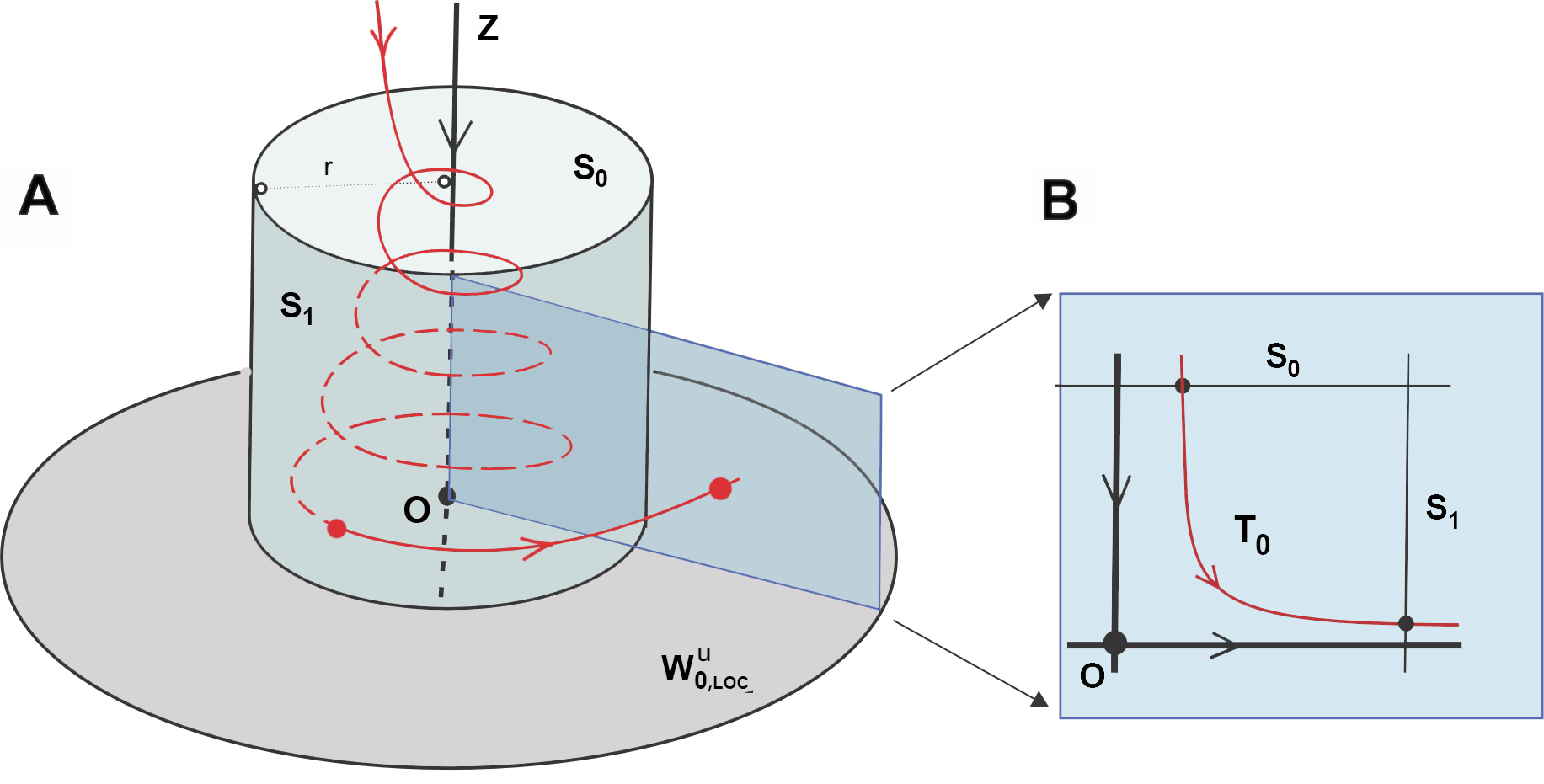}
\end{center}
\caption{(A) A half-neighborhood of the saddle-focus $O$ is bounded by its local unstable manifold $z=0$ and the cross-sections $S_0: \{z=0, r\leq 1\}$ and $S_1: \{z=1, r\leq 1\}$. (B) Orbits of the linear system~(\ref{eq2}) define the return map $T_0: S_0 \to S_1$.}  
\label{fig4}
\end{figure}

Consider a half-neighborhood of $O$, as in Fig.~\ref{fig4}A. Take the cross-sections: the disk $S_0: \{ z=1,\, r\le 1 \}$, and the cylinder $S_1: \{ r=1,\, 0\le z\le 1 \}$. One can construct a map $T_0: S_0 \to S_1$, see Fig.~\ref{fig4}B. The solution of 
(\ref{eq3}) is 
$$
 r  = r_0 e^{\rho t} , \qquad \varphi = \omega t + \varphi_0, \qquad  z= e^{-\lambda\,t}. $$
The transition time $\tau$ between the cross-sections is found from 
$$1=r_0 e^{\rho \tau},$$
hence 
$$\tau=-\frac{1}{\rho} \ln (r_0),$$
and the map $T_0$ is given by 
$$\begin{aligned}
z_1 & = r_0^{\nu},\\
\varphi_1 & = \frac{\omega}{\rho}\ln \frac{1}{r_0} + \varphi_0, \quad  \mbox{mod}~ 2\pi,
\end{aligned}$$ 
[where $(r_0,\varphi_0)$ are the coordinates of the initial point in $S_0$ and $(z_1,\varphi_1)$ are the coordinates of its
image by $T_0$ in $S_1$, i.e., the point where the trajectory of $(r_0,\varphi_0)$ hits the cross-section $S_1$]. 
The main condition is 
$$\nu = \frac{\lambda}{\rho} >1.$$

Consider the global map $T_1$
$$\begin{aligned}
r_0  & = f(z_1, \varphi_1, \mu), \\
\varphi_0 & = g(z_1, \varphi_1, \mu), \quad  \mbox{mod}~2\pi.
\end{aligned}$$
[That is, we assume that all orbits of the nonlinear system (\ref{eq1}), which leave the neighborhood of the saddle-focus $O$, return to the cross-section $S_0$. The map $T_1$ relates the exit point in $S_1$ with the point where the orbit lands to $S_0$. The unstable manifold $W^u_O$ also gets to $S_0$ --
the region bounded by the piece of $W^u_O$ before it lands to $S_0$ and the disc in $S_0$ that is bounded by the circle $T_1(W^u_{O,loc}\cap S_1) \cap S_0$ is the Shilnikov funnel. The outside orbits which enter the funnel (through $S_0$) never leave it. The parameter $\mu$ is introduced to allow one to study bifurcations in the funnel.]

\begin{figure}[t!]
\begin{center} 
\includegraphics[width=0.45\textwidth]{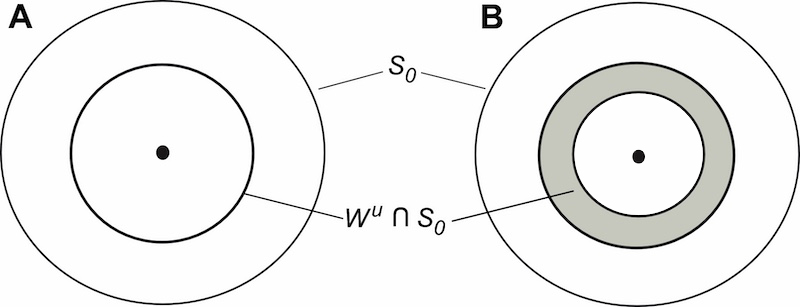}
\end{center}
\caption{At $n=1$, the global map sends (A) $W^u_{O,loc}\cap S_1$ to a circle surrounding the origin in the cross-section $S_0$, and (B) the entire cross-section $S_1$ to an annulus surrounding the origin.}  
\label{fig5}
\end{figure}

The image $W^u_O\cap S_0$ will be found from here 
\begin{equation}\label{eq7}
\begin{aligned}
r_0 & = f(0, \varphi_1, \mu), \\
\varphi_0 & = g(0, \varphi_1, \mu)=n \varphi_1+g_1(0,\varphi_1, \mu)  \; \mbox{mod} ~2\pi, 
\end{aligned}
\end{equation}
where $n$ is an integer, and $f$ and $g_1$ are periodic 
in $\varphi_1$. 
In the case of a system in $R^3$, the integer $n$ can be either $0$, or $+1$. 

If $f(0, \varphi^*, \mu)$ vanishes for some $\varphi^*$, then the system has a homoclinic loop to the saddle-focus, see Fig.~\ref{fig3}. When $\nu>1$ this implies the presence of complex dynamics, at this value of $\mu$ and for all close values of $|\mu|$ [by Shilnikov theorem \cite{LP1,LP2,LP3}]. We therefore impose that
$$f(0, \varphi, \mu) >0.$$
[Our primary goal is to address the existence of a torus as an attractor, which requires the absence of a homoclinic loop. On the other hand, we want to remain close to the homoclinic regime. For this reason, the model maps introduced below assume that $f_1(0, \varphi, \mu)$ is small.]

~\\
Now, two cases are possible: \\
$W_O^u \cap S_0$ encloses the point $(0,0)$; see Fig.~\ref{fig5}A. In this case, $n=1$ in Eq.~(\ref{eq7}). Therefore, the image of $S_1$ is an annulus around $(0,0)$, with the outer and inner boundaries given by the images of $z_1=1$ and $z_1=0$, respectively. \\
$W_O^u \cap S_0$ does not enclose $(0,0)$; see Fig.~\ref{fig6}A. In this case, $n=0$, and the image of $S_1$ is a ring, as shown in Fig.~\ref{fig6}B. Here, it is necessary to require that $(0,0) \not\in T_1 S_1$.

{\em Case $n=1$}. Let us assume that the the global map $T_1$ has the following form: 
\begin{figure}[t!]
\begin{center} 
\includegraphics[width=0.45\textwidth]{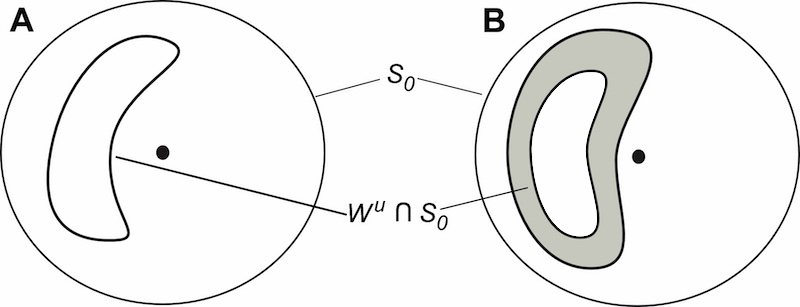}
\end{center}
\caption{At $n=0$, the global map sends (A) $W^u_{O,loc}\cap S_1$ to a circle that does not enclose the origin, and (B) the entire cross-section $S_1$ to an annulus that does not contain the origin.}  
\label{fig6}
\end{figure}
$$
\begin{aligned}
r_0  & = \mu \alpha (\varphi_1) + O(z), \\
\varphi_0 & = \varphi_1 + \varphi^* + \mu \beta (\varphi_1) + O(z)  \; \mbox{mod}~2\pi, 
\end{aligned}
$$
where $\alpha(\varphi_1)>0$. 
[The main idea behind this formula is that the map $T_1$ degenerates at $\mu=0$: the image of the circle $W^u_{O,loc}:, z=0$ collapses to a point ($r=0$). This assumption facilitates the analysis of dynamics recurrent to a neighborhood of the saddle-focus for small $\mu>0$. It becomes feasible in several examples of homoclinic configurations, as discussed at the end of the text for systems with rotational symmetry.]

The combined map 
$T=T_0\circ T_1$ will be given by 
$$\begin{aligned}
\bar z  & =  \left ( \mu \alpha (\varphi) + O(z)\right )^\nu, \\
\bar \varphi & = \varphi + \frac{\omega}{\rho} 
\ln \frac{1}{\mu \alpha (\varphi)+O(z)} + \varphi^*+ \mu \beta (\varphi)+ O(z),\ \mbox{mod}~2\pi.
\end{aligned}$$
[The map $T$ takes a point $(z,\varphi)\in S_1$ to
$(\bar z, \bar \varphi) \in S_1$. This map, by construction, completely determines the dynamics inside the funnel.]

Let us make the coordinate transformation $z\to \mu^\nu z$. Then the map $T$ recasts as 
\begin{equation}
\begin{aligned}
\bar z  & =  \alpha(\varphi)^\nu+\cdots, \\
\bar \varphi & = \varphi + \frac{\omega}{\rho} \ln \frac{1}{\alpha (\varphi)}+ \widetilde{\omega}(\mu) +\cdots    \;\; \mbox{mod}~2\pi, 
\end{aligned}
\label{eq10}
\end{equation}
where $\widetilde{\omega}(\mu) \to \infty$ as $\mu \to 0$, and the dots stand for the terms going to zero as $\mu \to 0$ along with their first derivatives [this follows from the condition $\nu>1$].

One can see that the circle map 
\begin{equation}
\bar \varphi  = \varphi + \frac{\omega}{\rho} \ln \frac{1}{\alpha (\varphi)}+ \widetilde{\omega}(\mu)  \; \mbox{mod}~2\pi, 
\label{eq11}
\end{equation}
is a diffeomorphism when 
  \begin{equation}
\frac{\omega}{\rho} \frac{\alpha^\prime(\varphi)}{\alpha(\varphi)} <1. 
\label{eq12}
\end{equation}

{\bf Proposition~I.} {\em If the inequality~(\ref{eq12}) holds, then for all sufficiently small $\mu>0$ the map~(\ref{eq10}) has a stable, invariant, and smooth closed curve of the form
$z=h(\varphi;\mu)$.}\\

Thus, under quite broad conditions, system~(\ref{eq1}) has a stable, invariant, smooth two-dimensional torus. However, complex dynamics is also possible.\\

{\bf Proposition II.}\, {\em Let there exist an interval $I= [ \varphi_1,\, \varphi_2]$ such that, for some $m\ge 2$, either
$$\alpha^\prime(\varphi)<0 \qquad \mbox{for all}\;\; \varphi\in I
\qquad \mbox{and}$$
$$
\frac{\omega}{\rho} \ln \frac{\alpha(\varphi_1)}{\alpha(\varphi_2)} > 2 \pi(m+1), 
$$
or 
$$
\frac{\omega}{\rho} \frac{\alpha^\prime(\varphi)}{\alpha(\varphi)} > 2, \qquad \mbox{for all}\;\; \varphi\in I
\qquad \mbox{and}$$
$$
\frac{\omega}{\rho} \ln \frac{\alpha(\varphi_2)}{\alpha(\varphi_1)} > 2 (\varphi_2 - \varphi_1)  + 2\pi(m+1).
$$
Then, for all small enough $\mu>0$, the map $T$ possesses a hyperbolic set that is topologically conjugate to the Bernoulli shift with $m$-symbols.}\\

\begin{figure}[t!]
\begin{center} 
\includegraphics[width=0.25\textwidth]{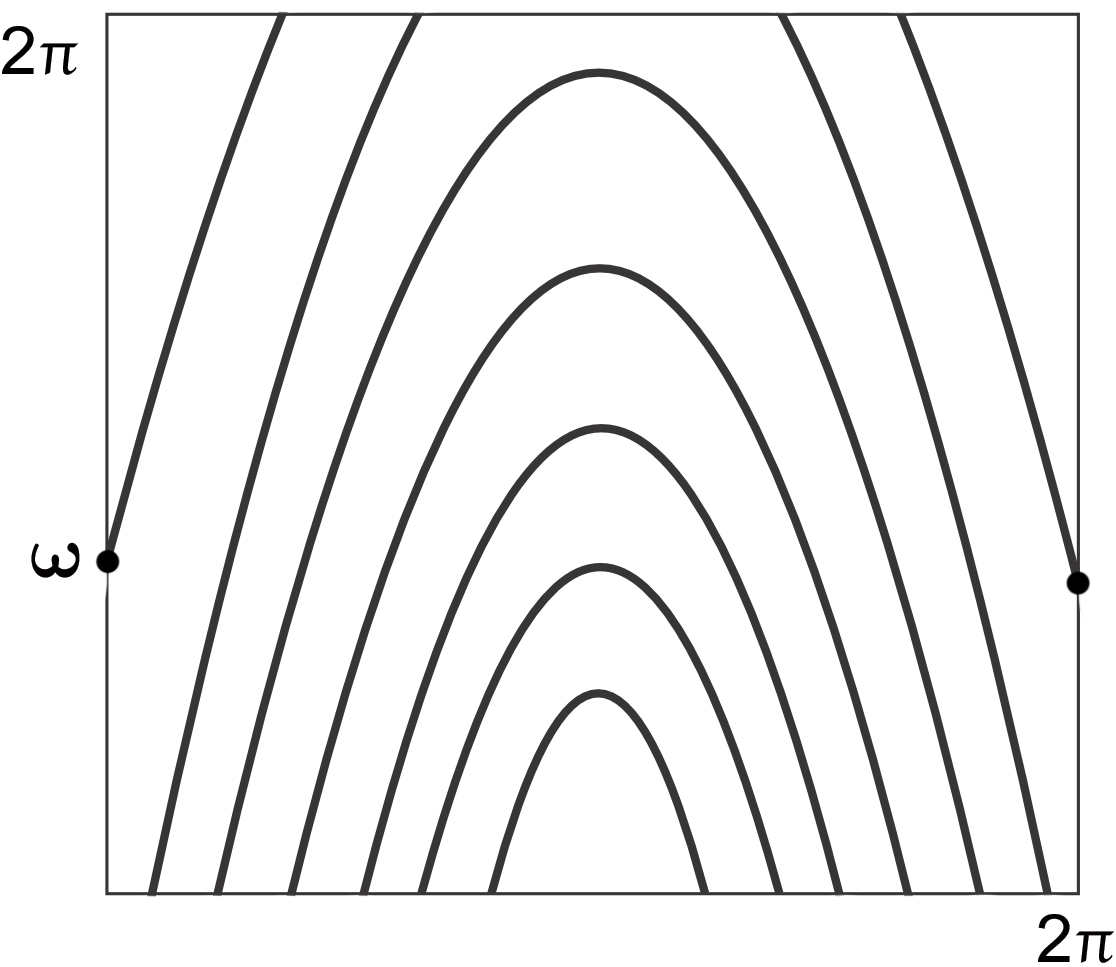}
\end{center}
\caption{For large $A$, the map $\varphi\mapsto h(\varphi) = A\sin\varphi+\widetilde{\omega}(\mu)$ has many expanding branches, which creates non-trivial hyperbolic sets, inherited by map $T$ for small $\mu>0$, see Eq.~(\ref{eq20}).}  
\label{fig7}
\end{figure}

Note that maps of type (\ref{eq10}) were studied by V.S.~Afraimovich and L.P.~Shilnikov in connection with the problem of periodic forcing applied to an oscillatory system that is assumed to possess a stable periodic orbit emerging from a separatrix loop of a saddle equilibrium [see \cite{ALP3,ALP1,ALP2}, as well as \cite{sync} and references therein]. In the case where $\alpha(\varphi)=1+a,\sin(\varphi)$, the main resonant zones or Arnold tongues were studied in detail, together with the breakdown of an invariant torus accompanied by various bifurcation phenomena, including period-doubling cascades, homoclinic tangencies, and global bifurcations involving a saddle-node fixed point and a non-smooth closed invariant curve. Additional dynamical effects arising from the overlap of resonant zones were also described.
\\

{\em The case $n=0$.}
Assume that the global map $T_1$ has the form
$$\begin{aligned}
r_0 &  = \mu \alpha (\varphi_1) + O(z), \\
\varphi_0 & = \varphi^* + \mu \beta (\varphi_1)+\O(z) \; \mbox{mod}~2\pi. 
\end{aligned}
$$
Then $T = T_0 \circ T_1$ will be 
$$\begin{aligned}
\bar z & =  (\mu \alpha (\varphi) + O(z))^\nu, \\
\bar \varphi & = \frac{\omega}{\rho} \ln \frac{1}{\mu \alpha (\varphi)+O(z)} + \varphi^* + \mu \beta (\varphi)+ O(z)    \; \mbox{mod}~2\pi,\end{aligned} 
$$
where $\alpha(\varphi_1)>0$. 

\begin{figure}[t!]
\begin{center} 
\includegraphics[width=0.25\textwidth]{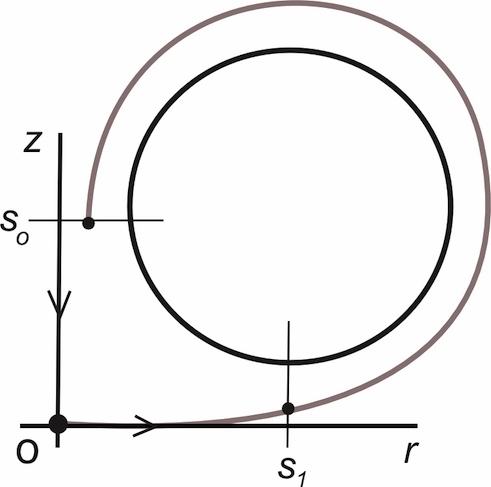}
\end{center}
\caption{A stable periodic orbit emerging from a homoclinic loop to a saddle with a negative saddle value.}  
\label{fig8}
\end{figure}

Make a coordinate transformation $z \to \mu^\nu z$. Then, we obtain 
\begin{equation}\label{eq20}
\begin{aligned}
\bar z &  = \alpha(\varphi)^\nu +\cdots, \\ 
\bar \varphi  & = \frac{\omega}{\rho} \ln \frac{1}{\alpha (\varphi)} + \widetilde{\omega}(\mu) + \cdots,    \;\; \mbox{mod}~2\pi,
\end{aligned}
\end{equation} 
where  $\widetilde{\omega}(\mu) \to \infty$ as $\mu \to 0$, and the dots stand for terms going to zero as $\mu \to 0$, along with their first derivatives.

The corresponding characteristic equation is given by 
$$\left | 
\begin{array}{cc}
  o(1) - \lambda  &~~~~~ \nu \alpha^{\nu-1}(\varphi)\alpha^\prime(\varphi) + o(1)  \\
  o(1) &  - \frac{\omega}{\rho} \frac{\alpha^\prime (\varphi)}{\alpha(\varphi)}+ o(1) -\lambda
\end{array} \right  |=0.
$$
One root is $o(1)$ at small $\mu$, while the other is 
$\lambda=-\frac{\omega}{\rho} \frac{\alpha^\prime (\varphi)}{\alpha(\varphi)}$. Therefore, under the condition
\begin{equation}\label{eq21}
\left|\frac{\omega}{\rho} \frac{\alpha^\prime (\varphi)}{\alpha(\varphi)}\right| < 1
\end{equation}
the fixed point is stable.\\
 
{\bf Proposition III.} {\em The system of differential equations (\ref{eq1}) has a stable periodic orbit.}\\

The map~(\ref{eq20}) is close to the one-dimensional map $\bar\varphi = f(\varphi)+\widetilde{\omega}(\mu)$ where we put $f(\varphi)=\frac{\omega}{\rho} \ln \frac{1}{\alpha(\varphi)}$.
Letting $f(\varphi)=A \sin(\varphi)$, for example, we obtain complex dynamics for sufficiently large values of $A$, see Fig.~\ref{fig7}. Similarly to the proof of Proposition II, one establishes that the original 2D map (\ref{eq20}) at such $A$ possesses a hyperbolic set conjugate to a Bernoulli shift.\\

 \begin{figure}[t!]
\begin{center} 
\includegraphics[width=0.29\textwidth]{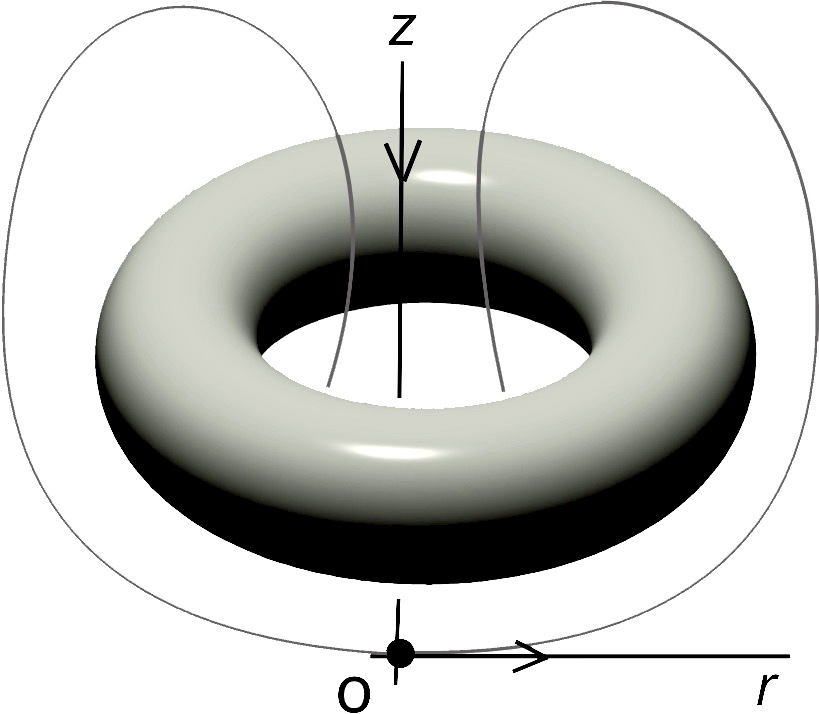}
\end{center}
\caption{A 2D torus in a system with the rotational symmetry around the $z$-axis.}  
\label{fig9}
\end{figure}

\centerline{\em Geometric interpretations}

Let us consider a system in a plane whose phase portrait is as sketched in Fig.~\ref{fig8}. Here, $O$ is a saddle with a negative saddle value [which is the sum of its characteristic exponents]. Then, the Poincar\'e map can have the following form    
$$
\bar z =   a z^\nu + \alpha \mu,
$$
where $\nu>1$, and $\mu$ is a small parameter. This map has a stable fixed point at $\mu>0$. Let us append our 2D system with a third equation $\dot\varphi= \omega$. We obtain a 3D system with a rotational symmetry, which has a two-dimensional stable torus, see Fig.~\ref{fig9}. This torus will also exist at small $\mu>0$ for small perturbations [that may break the rotational symmetry] such that $\alpha(\varphi)>0$.  

\begin{figure}[b!]
\begin{center} 
\includegraphics[width=0.42\textwidth]{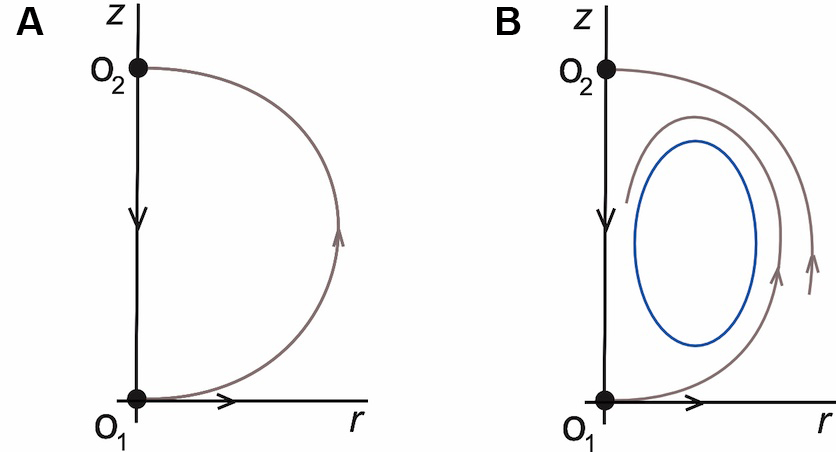}
\end{center}
\caption{Breakdown of a heteroclinic cycle between two saddles (A) leading to the emergence of a stable periodic orbit (blue orbit in panel~B). }  
\label{fig10}
\end{figure}

Another example of a two-dimensional system is sketched in Fig.~\ref{fig10}A. Here $O_1$ has the saddlei index $\nu_1>1$ and $O_2$ has $\nu_2>1$. The small splitting of the heteroclinic orbit $\Gamma_1$ as shown in  Fig.~\ref{fig10}B will also create a stable periodic orbit. As before, the 3D system with $\dot\varphi= \omega$ will have a 2D stable torus. It will also persist under sufficiently small perturbations.     

\begin{figure}[b!]
\begin{center} 
\includegraphics[width=0.25\textwidth]{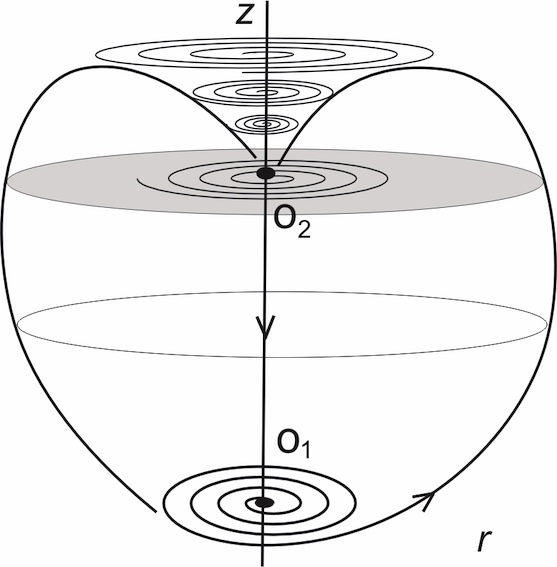}
\end{center}
\caption{A heteroclinic cycle between the saddle-focus $O_1$ and the saddle-node $O_2$: all orbits in the 2D unstable manifold $W^u_{O_1}$ tend to the saddle-node as $t\to+\infty$, and the 1D unstable manifold of $O_2$ merges with the 1D stable manifold of $O_1$.}  
\label{fig11}
\end{figure}

~\\
{\em The second interpretation}. Consider a three-dimensional system with a saddle-focus $O_1$ that has a one-dimensional stable manifold, and a saddle-node $O_2$ that attracts all trajectories from $W^u_{O_1}$. At the critical moment, the one-dimensional unstable manifold of the saddle-node enters the saddle-focus $O_1$ as $t \to +\infty$; see Fig.~\ref{fig11}. This configuration represents a global bifurcation of codimension three: one parameter controls the saddle-node, while two additional parameters govern the merging of $W^u_{O_2}$ with $W^s_{O_1}$. If the system possesses a rotational symmetry, only one parameter is required to create the saddle-node.
By introducing small perturbations that destroy the saddle-node, one arrives at the class of models described in this paper. In this case, the global map $T_1$ acts as illustrated in Fig.~\ref{fig5}B or in Fig.~\ref{fig6}B.

\begin{acknowledgements}
A.L.S. is grateful to Dima Turaev for his always insightful suggestions. He also thanks the members of his NeurDS Lab:  C. Hinsley, Y. Keleta and Dr. J. Scully for always fruitful and helpful discussions during Friday lab meetings.  
\end{acknowledgements}


\end{document}